\documentclass[final,3p,times]{elsarticle}
\usepackage{amsmath,amssymb,amscd,amsthm,verbatim,alltt,amsfonts,array}
\usepackage[english]{babel}
\usepackage{latexsym}
\usepackage{amssymb}
\usepackage{euscript}
\usepackage{graphicx}
\usepackage{pgf,tikz}
\usetikzlibrary{arrows}
\usepackage{mathrsfs}
\definecolor{qqffqq}{rgb}{0.,1.,0.}
\definecolor{qqqqff}{rgb}{0.,0.,1.}
\definecolor{zzwwqq}{rgb}{0.6,0.4,0.}
\definecolor{ffffqq}{rgb}{1.,1.,0.}
\definecolor{ttffqq}{rgb}{0.2,1.,0.}
\definecolor{eqeqeq}{rgb}{0.8784313725490196,0.8784313725490196,0.8784313725490196}
\definecolor{qqzzff}{rgb}{0.,0.6,1.}
\definecolor{ffqqqq}{rgb}{1.,0.,0.}
\definecolor{cqcqcq}{rgb}{0.7529411764705882,0.7529411764705882,0.7529411764705882}
\definecolor{zzttqq}{rgb}{0.6,0.2,0.}
\definecolor{qqqqff}{rgb}{0.,0.,1.}
\definecolor{xdxdff}{rgb}{0.49019607843137253,0.49019607843137253,1.}
\definecolor{qqffqq}{rgb}{0.,1.,0.}
\definecolor{ffffff}{rgb}{1.,1.,1.}
\definecolor{ffqqqq}{rgb}{1.,0.,0.}
\definecolor{ccffcc}{rgb}{0.8,1.,0.8}
\definecolor{ffffqq}{rgb}{1.,1.,0.}
\definecolor{qqffqq}{rgb}{0.,1.,0.}
\definecolor{uuuuuu}{rgb}{0.26666666666666666,0.26666666666666666,0.26666666666666666}
\definecolor{ccqqww}{rgb}{0.8,0.,0.4}
\definecolor{bfffqq}{rgb}{0.7490196078431373,1.,0.}
\definecolor{ffffff}{rgb}{1.,1.,1.}
\definecolor{ttqqqq}{rgb}{0.2,0.,0.}
\definecolor{qqffff}{rgb}{0.,1.,1.}
\usepackage{amsmath,amssymb,amsfonts,epsf,amsthm}
\usepackage{hyperref}
\hypersetup{colorlinks=true, urlcolor=blue, citecolor=cyan, pdfborder={0 0 0},}
\usepackage{soul}
\usepackage{natbib}
\usepackage{multirow}

\usepackage[ruled, vlined, Algorithm]{algorithm}
\usepackage{algorithmic}

\newtheorem{theorem}{Theorem}[section]
\newtheorem{lemma}{Lemma}[section]
\newtheorem{proposition}{Proposition}[section]

\newtheorem{definition}{Definition}[section]

\numberwithin{equation}{section}
\usepackage{amssymb}
\journal{--}
\allowdisplaybreaks[1]
\begin{document}
\begin{frontmatter}
\title{ Some resolving parameters in a class of Cayley graphs }
\author[label1]{Jia-Bao Liu}
\ead{liujiabaoad@163.com;liujiabao@ahjzu.edu.cn}
\author[label2]{Ali Zafari \corref{1}}
\ead{zafari.math.pu@gmail.com; zafari.math@pnu.ac.ir}
\address[label1]{School of Mathematics and Physics,
Anhui Jianzhu University, Hefei 230601, P.R. China}
\address[label2]{Department of Mathematics, Faculty of Science,
Payame Noor University, P.O. Box 19395-4697, Tehran, Iran}
\cortext[1]{Corresponding author}
\begin{abstract}
Resolving parameters is a fundamental area of combinatorics with applications not only to many branches of combinatorics but also to other sciences. In this article, we construct a class of Toeplitz graphs, and will be denoted by $T_{2n}(W)$, so that they are Cayley graphs. First, we review some of the features of this class of graphs. In fact, this class of graphs are vertex transitive, and  by calculating the spectrum of the adjacency matrix related with them, we show that this class of  graphs cannot be edge transitive. Moreover, we show that  this class of graphs cannot be distance regular, and since the computing resolving parameters of a class of graphs such that are not distance regular is more difficult, then we regard this as justification for our focus on some resolving parameters. In particular, we determine the minimal resolving set, doubly resolving set and strong metric dimension for this class of graphs.
\end{abstract}
\end{frontmatter}
\section{Introduction}
\label{sec:introduction}
The graphs in this paper are simple, undirected and connected. An automorphism of a graph $\Gamma$ is a permutation $\varphi$ of the vertex set of $\Gamma$ with the property that, for any vertices $x$ and $y$ we have $x$ is adjacent to $y$ in $\Gamma$ if and only if $\varphi(x)$ is adjacent to  $\varphi(y)$ in $\Gamma$. The set of all automorphisms of a graph $\Gamma$, with the operation of composition of permutations, is a permutation group on $V(\Gamma)$, a subgroup of the symmetric group on $V(\Gamma)$. This is the automorphism group of $\Gamma$, denoted by $Aut(\Gamma)$. Suppose $\Gamma_1$ and $\Gamma_2$ are two graphs. If there is a bijection, $\varphi$ say, from $V(\Gamma_1)$ to $V(\Gamma_2)$ so that $x$ is adjacent to $y$ in $\Gamma_1$ if and only if $\varphi(x)$ is adjacent to  $\varphi(y)$ in $\Gamma_2$, then we say that $\Gamma_1$ is isomorphic to $\Gamma_2$.
If we consider a graph $\Gamma$ as a network, then the network stability is very important to us, and especially, if  graph $\Gamma$ is vertex transitive, that is, $Aut(\Gamma)$ acts transitively on $V(\Gamma)$, then the cost of study the network will be very low, and hence the  network will be more stable.
Consider a finite group $G$, and suppose $Q$ is a subset  of $G$ so that it is closed under taking inverses and does not contain the identity, then the Cayley graph $\Gamma=Cay(G, Q)$ has vertex set $G$  and edge set
 $E(\Gamma) = \{\{x, y\} \, | \,\, x^{-1}y \in Q\}$.
Thus the studying of Cayley graphs is very useful, because every Cayley graph is vertex transitive [1]. 
The distance between any pair $u, v \in V(\Gamma)$ of vertices of $\Gamma$ is the length of geodesic between $u$ and $v$, denoted by $d_{\Gamma}(u, v)$ or simply $d(u, v)$. A vertex $x\in V(\Gamma)$ is said to resolve a pair $u, v \in V(\Gamma)$ if $d(u, x)\neq d(v, x)$.
Resolving parameters is a fundamental area of combinatorics with applications not only to many branches of combinatorics but also to other sciences. For an arranged subset $R = \{r_1, r_2, ..., r_m\}$ of vertices in a connected graph $\Gamma$  the metric representation of a vertex $v$ in $\Gamma$, is the $m$-vector $r(v | R) = (d(v, r_1), d(v, r_2), ..., d(v, r_m ))$  relative to $R$. Also, the subset $R$ is considered as resolving set for $\Gamma$ if any pair of vertices of $\Gamma$ is distinguished by some vertices of $R$. A resolving set with least number of vertices is referred as metric basis for $\Gamma$ and the cardinality of such resolving set is know as metric dimension denoted by $\beta(\Gamma)$. The metric dimension of a graph $\Gamma$ is the least number of vertices in a set with the property that the list of distances from any vertex to those in the set uniquely identifies that vertex. The concept of  the metric dimension in algebraic graph theory date back to the 1970s. It was defined  independently by Harary and Melter [2] 
 and by Slater [3].
In recent years, a considerable literature has developed [4]. 
This concept has different applications in the areas of network discovery and verification [5]. 
For more details see [6-9].
An (n x n) matrix $T = (t_{ij})$ is called a Toeplitz matrix if $t_{ij} = t_{i+1,j+1}$ for each $i,j = 1, ..., n-1$, see [10]. 
In fact, a Toeplitz matrix is a square matrix so that entries in every diagonal parallel to the main diagonal are equal, and hence a Toeplitz matrix is determined by its first row and column. A simple undirected graph $\Gamma$ with vertex set $\{1, ..., n\}$ and its adjacency matrix $T= (t_{ij})$ is called a Toeplitz graph if  $T$ is the  Toeplitz matrix. In this paper, we consider a class of Toeplitz graphs, will be denoted by  $T_{2n}(W)$, so that  they  are  Cayley  graphs as follows:\newline
Let $n$ be a fixed  even integer is greater than or equal $4$, also, let $[2n]= \{1, 2, ..., 2n\}$ and $[x_{2n}]= \{x_1, x_2, ..., x_{2n}\}$ be  corresponding sets so that $x_i=i$. Hence we say that $x_i<x_j$, if $i<j$. Now, let $W_1=\{x_1, x_3, ..., x_{2n-1}\}$ and $W_2=\{x_n \}$ be  subsets of the set $[x_{2n}]$, and let $W=W_1\cup W_2=\{x_1, x_3, ..., x_n, ..., x_{2n-1}\}$ be a refinement of union of the two  sets $W_1$ and $W_2$ so that
$1=x_1< x_3< ...<x_n<...<x_{2n-1}$. We can see that   a graph with $2n$ vertices so that the vertices are labelled  by the set $\{1, 2, ..., 2n\}$, and  the edge set
$\{ij |\,\ i, j \in [2n],  |j-i|=x_t \,\ \text{for some $x_t\in W$ }\}$ is Toeplitz graph $T_{2n}(W)$. For more result of the Toeplitz graphs see [11, 12]. 
Figure 1 shows the Toeplitz graph $T_{8}(1, 3, 4, 5, 7)$.
\begin{center}
\begin{tikzpicture}[line cap=round,line join=round,>=triangle 45,x=5.0cm,y=5.0cm]
\clip(2.97703675890742,1.8818489102405913) rectangle (5.48467447718966,4.0084674699145735);
\fill[color=ffffff,fill=ffffff,fill opacity=0.1] (4.,2.5) -- (4.5,2.5) -- (4.853553390593274,2.853553390593274) -- (4.853553390593274,3.353553390593274) -- (4.5,3.707106781186548) -- (4.,3.707106781186548) -- (3.646446609406726,3.353553390593274) -- (3.646446609406726,2.853553390593274) -- cycle;
\draw [color=ffffff] (4.,2.5)-- (4.5,2.5);
\draw [color=ffffff] (4.5,2.5)-- (4.853553390593274,2.853553390593274);
\draw [color=ffffff] (4.853553390593274,2.853553390593274)-- (4.853553390593274,3.353553390593274);
\draw [color=ffffff] (4.853553390593274,3.353553390593274)-- (4.5,3.707106781186548);
\draw [color=ffffff] (4.5,3.707106781186548)-- (4.,3.707106781186548);
\draw [color=ffffff] (4.,3.707106781186548)-- (3.646446609406726,3.353553390593274);
\draw [color=ffffff] (3.646446609406726,3.353553390593274)-- (3.646446609406726,2.853553390593274);
\draw [color=ffffff] (3.646446609406726,2.853553390593274)-- (4.,2.5);
\draw [color=ffffff] (4.,3.707106781186548)-- (4.5,3.707106781186548);
\draw [color=ffffff] (4.,3.707106781186548)-- (4.5,3.707106781186548);
\draw [color=ffffff] (4.5,3.707106781186548)-- (4.853553390593274,3.353553390593274);
\draw (4.,3.707106781186548)-- (4.5,3.707106781186548);
\draw (4.5,3.707106781186548)-- (4.853553390593274,3.353553390593274);
\draw (4.853553390593274,3.353553390593274)-- (4.853553390593274,2.853553390593274);
\draw (4.853553390593274,2.853553390593274)-- (4.5,2.5);
\draw (4.5,2.5)-- (4.,2.5);
\draw (4.,2.5)-- (3.646446609406726,2.853553390593274);
\draw (3.646446609406726,2.853553390593274)-- (3.646446609406726,3.353553390593274);
\draw (3.646446609406726,3.353553390593274)-- (4.,3.707106781186548);
\draw (4.,3.707106781186548)-- (4.853553390593274,2.853553390593274);
\draw (4.,3.707106781186548)-- (4.,2.5);
\draw (4.5,3.707106781186548)-- (4.5,2.5);
\draw (4.5,3.707106781186548)-- (3.646446609406726,2.853553390593274);
\draw (4.853553390593274,3.353553390593274)-- (4.,2.5);
\draw (4.853553390593274,3.353553390593274)-- (3.646446609406726,3.353553390593274);
\draw (4.853553390593274,2.853553390593274)-- (3.646446609406726,2.853553390593274);
\draw (4.5,2.5)-- (3.646446609406726,3.353553390593274);
\draw [color=qqqqff] (4.,3.707106781186548)-- (4.5,2.5);
\draw [color=qqqqff] (4.5,3.707106781186548)-- (4.,2.5);
\draw [color=qqqqff] (4.853553390593274,3.353553390593274)-- (3.646446609406726,2.853553390593274);
\draw [color=qqqqff] (4.853553390593274,2.853553390593274)-- (3.646446609406726,3.353553390593274);
\draw (3.605789681405645,2.283371676873298) node[anchor=north west] {Figure 1. The Toeplitz graph $T_{8}(1, 3, 4, 5, 7)$};
\begin{scriptsize}
\draw [fill=black] (4.,2.5) circle (1.5pt);
\draw[color=black] (3.9871805747525624,2.4557817288284704) node {$1$};
\draw [fill=black] (4.5,2.5) circle (1.5pt);
\draw[color=black] (4.50997430400575,2.4546426394937788) node {$2$};
\draw [fill=uuuuuu] (4.853553390593274,2.853553390593274) circle (1.5pt);
\draw[color=uuuuuu] (4.903019837271164,2.900853636751041) node {$3$};
\draw [fill=uuuuuu] (4.853553390593274,3.353553390593274) circle (1.5pt);
\draw[color=uuuuuu] (4.903019837271164,3.388203723342995) node {$4$};
\draw [fill=uuuuuu] (4.5,3.707106781186548) circle (1.5pt);
\draw[color=uuuuuu] (4.50997430400575,3.7589710779390247) node {$5$};
\draw [fill=uuuuuu] (4.,3.707106781186548) circle (1.5pt);
\draw[color=uuuuuu] (4.022624217413796,3.7589710779390247) node {$6$};
\draw [fill=uuuuuu] (3.646446609406726,3.353553390593274) circle (1.5pt);
\draw[color=uuuuuu] (3.5964132201565324,3.3793428126776868) node {$7$};
\draw [fill=uuuuuu] (3.646446609406726,2.853553390593274) circle (1.5pt);
\draw[color=uuuuuu] (3.5964132201565324,2.8919927260857325) node {$8$};
\end{scriptsize}
\end{tikzpicture}
\end{center}
In particular, we can verify that the Toeplitz graph $T_{2n}(W)$ which is defined already is isomorphic to the Cayley graph $\Lambda=Cay(\mathbb{D}_{2n}, \Psi)$, where $$\mathbb{D}_{2n}=<a,b \,\  | \,\, a^n=b^2=1  , ba=a^{n-1}b>,$$ is the dihedral group of order $2n$, and
$\Psi=\{ab, a^{2}b, ... , a^{n-1}b, b\}\cup\{a^{\frac{n}{2}}\}$ is an inverse closed subset of $\mathbb{D}_{2n}-\{1\}$. Thus, the Toeplitz graph $T_{2n}(W)$ is a vertex transitive.
Also for convenience, we can use the symbols in the Cayley graph $\Lambda=Cay(\mathbb{D}_{2n}, \Psi)$, instead of the symbols in the Toeplitz graph $T_{2n}(W)$. Some metrics for  a class of distance regular graphs computed in [13,14]. 
On the other hand, computing  metrics  of a class of graphs such that are not distance regular is more difficult, and hence we regard this as justification for our focus on some resolving parameters in the Cayley graph $\Lambda=Cay(\mathbb{D}_{2n}, \Psi)$. The important results of this article will be presented in two sections 3.1 and 3.2. In section 3.1, first, we will be determining the automorphism group of the Cayley graph $\Lambda=Cay(\mathbb{D}_{2n}, \Psi)$, also, we will show  that the Cayley graph $\Lambda=Cay(\mathbb{D}_{2n}, \Psi)$ cannot be  distance regular. In particular, we will prove that the Cayley graph $\Lambda=Cay(\mathbb{D}_{2n}, \Psi)$ cannot be edge transitive. Moreover, in section 3.2, we will be computing  some resolving parameters for this class of Cayley graphs.\newline
\section{Definitions And Preliminaries}
\begin{definition} (see [15]). \label{b.1} 
A graph $\Gamma$ is edge transitive if its automorphism group acts transitively on $E(\Gamma)$.
\end{definition}
\begin{definition} (see [15]). \label{b.2} 
A graph $\Gamma$ is 1-transitive or symmetric if its automorphism group acts transitively on the set of  paths of length 1 or 1-arcs.
\end{definition}
\begin{proposition} (see [15]). \label{b.3} 
Let $\Gamma$ be a symmetric  graph of valency $k$, and let $\lambda$ be a simple eigenvalue of $\Gamma$, then $\lambda=\pm k$.
\end{proposition}
\begin{definition}  (see [16]). \label{b.4} 
Suppose that $\Gamma$ is a regular graph of valency $k$ and for any two vertices $u$ and $v$ in $\Gamma$, if $d(u, v)=r$, then we have $|\Gamma_{r+1}(v)\cap\Gamma_{1}(u)|=b_r$, and $|\Gamma_{r-1}(v)\cap\Gamma_{1}(u)|=c_r$  $(0\leq r\leq d)$. Then we say that $\Gamma$ is a distance regular graph.
\end{definition}
\begin{proposition} (see [16]). \label{b.5} 
If $\Gamma$ is a distance regular graph with diameter $d$, then $\Gamma$ has exactly $d+ 1$ distinct eigenvalues.
\end{proposition}
\begin{definition}  (see [17]).  \label{b.6} 
Suppose $\Gamma$ is a  graph of order at least 2, vertices $x, y\in V(\Gamma)$  are said to doubly resolve vertices $u, v\in V(\Gamma)$ if
$d(u, x) - d(u, y) \neq d(v, x) - d(v, y)$. A  set $S = \{s_1, s_2, ..., s_l\}$ of vertices of $\Gamma$ is a doubly resolving set of $\Gamma$ if every
two distinct vertices of $\Gamma$ are doubly resolved by some two vertices of $S$. A doubly resolving set with minimum cardinality is called minimal doubly resolving set. This minimum cardinality is denoted by $\psi(\Gamma)$.
\end{definition}
\begin{definition} (see [18]). \label{b.7} 
Let $\Gamma$ be a graph.   A vertex $w$  of $\Gamma$ strongly resolves two vertices $u$ and $v$ of $\Gamma$ if $u$ belongs to a
shortest $v - w$ path or $v$ belongs to a shortest $u - w$ path. A set $S= \{s_1, s_2, ..., s_m\}$ of vertices of $\Gamma$ is a
strong resolving set of $\Gamma$ if every two distinct vertices of $\Gamma$ are strongly resolved by some vertex of $S$.
The strong metric dimension of a graph $\Gamma$ is  the cardinality of  smallest  strong resolving set of $\Gamma$ and denoted by $sdim(\Gamma)$.
\end{definition}
\section{Main results}
\noindent \textbf{ 3.1 \,\ Some of the features of the Cayley graph $Cay(\mathbb{D}_{2n}, \Psi)$}\\

In this section we review some of the features of the Cayley graph $Cay(\mathbb{D}_{2n}, \Psi)$.
It is well known that, the spectrum of a graph is the spectrum of the adjacency matrix related with it, that is, its set of eigenvalues together with their multiplicities. If all the eigenvalues of the adjacency matrix of a graph  are integers, in this case, the graph related with it is called an integral graph,  see [19]. 
As we shall see, the theory of integral graphs has connections to some parts of graph theory, edge transitivity, and symmetric graph. In the next Theorem, we obtain the  automorphism group of the Cayley graph $Cay(\mathbb{D}_{2n}, \Psi)$ by applications of wreath product in graph theory, for more details of wreath product, see [20]. 
\begin{proposition}\label{c.1}
Let $n$ be an even integer greater than or equal $4$, and $\Lambda=Cay(\mathbb{D}_{2n}, \Psi)$ be a Cayley graph on the dihedral group $\mathbb{D}_{2n}$, where $\Psi$  which is defined already. If $k=\frac{n}{2}-1$, then
$Aut(\Lambda)\cong\mathbb{Z}_2 wr_{I} Sym{(k+1)}wr_{J} Sym{(2)}$, where $I=\{1, ... , k+1\}$ and $J=\{1, 2\}.$
\begin{proof}
We can see that the complement of $\Lambda$, denoted by   $\overline{\Lambda}$,  is isomorphic to the disjoint union of $2$ copies  of cocktail party graph $CP(\frac{n}{2} )$, and we can show that $CP(\frac{n}{2} )$ is isomorphic to the $Cay(\mathbb{Z}_n, S_k)$, where $\mathbb{Z}_n$ is the cyclic group of order $n$ and
$S_k=\{1, n-1, 2, n-2, ..., k, n-k\}$, see  [21, Proposition 3.2]. 
Hence given by  above discussion and the Theorem in [22], 
we have $Aut(\overline{\Lambda})\cong Aut(CP(\frac{n}{2} )) wr_{I} Sym{(2)}=\mathbb{Z}_2 wr_{I} Sym{( k+1))}wr_{J} Sym{(2)}$. In particular, we have $Aut(\overline{\Lambda})=Aut(\Lambda)$, because a simple undirected
graph and its complement have the same automorphism group.
\end{proof}
\end{proposition}
\begin{proposition}\label{c.2}
Let $n$ be an even integer greater than or equal $4$, and $\Lambda=Cay(\mathbb{D}_{2n}, \Psi)$ be a Cayley graph on the dihedral group $\mathbb{D}_{2n}$, where $\Psi$  which is defined already, then  $\Lambda$ cannot be a distance regular graph.
\end{proposition}
\begin{proof}
It is not hard to see that, the diameter of $\Lambda$ is $2$ and $\Lambda$ is not a bipartite graph, because $a^{\frac{n}{2}}\in \Psi$. Now by a similar way, which is done in proof of Proposition 11,
in [23] 
we can show that the adjacency matrix spectrum of $\Lambda$ is
$n+1, 1-n, 1^{(n-2)}, -1^{(n)}$, where the superscripts give the multiplicities of eigenvalues with multiplicity greater than one. Hence, $\Lambda$ has exactly four distinct eigenvalues. Moreover, based on Proposition 2.2, 
we know that if $\Lambda$ is a distance regular graph with diameter $d$, then $\Lambda$ has exactly $d+ 1$ distinct eigenvalues. Thus, $\Lambda$ cannot be a distance regular graph.
\end{proof}
\begin{proposition}\label{c.3}
Let $n$ be an even integer greater than or equal $4$, and $\Lambda=Cay(\mathbb{D}_{2n}, \Psi)$ be a Cayley graph on the dihedral group $\mathbb{D}_{2n}$, where $\Psi$  which is defined already, then  $\Lambda$ cannot  be a edge transitive graph.
\end{proposition}
\begin{proof}
By contradiction, suppose $\Lambda$ is a edge transitive graph. It is well known that a connected graph that is edge transitive and vertex transitive need not be 1-transitive. In particular, in  [24, p.59, 7.53],
Tutte proved that if a connected graph, regular of odd valency, is both vertex and edge transitive, then it is 1-transitive. Thus, if $\Lambda$ is a edge transitive graph, then it must be $\Lambda$ is a 1-transitive graph, because it is vertex transitive of odd valency $n+1$. On the other hand, based on Proposition 2.1, 
if $\lambda$  is a simple eigenvalue of a  1-transitive graph $\Lambda$, then $\lambda=\pm(n+1)$, which is not the case, see previous Proposition. This contradiction shows that $\Lambda$ cannot  be a edge transitive graph.
\end{proof}
\noindent \textbf{3.2   \,\ Metric dimension, minimal doubly resolving set, and strong resolving  set  of the Cayley graph $Cay(\mathbb{D}_{2n}, \Psi)$ }\\

\begin{theorem}\label{d.1}
If $n$ is an even integer greater than or equal $4$, and $\Lambda=Cay(\mathbb{D}_{2n}, \Psi)$ is a Cayley graph on the dihedral group $\mathbb{D}_{2n}$, where $\Psi$  which is defined already, then  the metric dimension of $\Lambda$ is $n$.
\end{theorem}
\begin{proof}
Let $V(\Lambda)=V_1 \cup V_2$, where $V_1=\{a,a^2, ..., a^{n}\}$ and $V_2=\{ab,a^2b, ..., a^nb\}$. For every pair of distinct vertices $x, y \in V (\Lambda)$, the length of a shortest path from $x$ to $y$ is   $d(x, y)=1$ or $2$, because the diameter of $\Lambda$ is $2$. In particular, if $R$ is an arranged subset of $V_1$ or $V_2$ in graph $\Lambda$ such that $|R|\leq n$, then we can show that $R$ is not a resolving set of $\Lambda$. Let $R=R_1\cup R_2$, be an arranged subset  of vertices in  graph $\Lambda$ such that  $R_1$ is a subset of $V_1$,  $R_2$  is a subset of $V_2$  and  $|R_1\cup R_2|= n$. In the following cases, we can be concluded that the metric dimension of $\Lambda$ is $n$.\newline

Case 1.
If $|R_1|\neq |R_2|$, then we can assume  that $|R_1|< |R_2|$. Hence, there is a pair of distinct vertices
$u_1, u_2 \in V(\Lambda)- R$,  such that $u_1, u_2\in V_1- R_1$, and a shortest path from $u_1$ to $u_2$ is
$d(u_1, u_2)=1$. Therefore, the metric representation of  the vertices $u_1, u_2\in V(\Lambda)-R$ is the same as $n$-vector, relative to $R$. Thus, $R$ is not a resolving set of $\Lambda$.
\newline

Case 2.
If $|R_1|=|R_2|$ and there are vertices  $x, y \in R_1$,  such that $x$ is adjacent to  $y$ in $\Lambda$, then there are vertices $u, v \in V_1-R_1$ such that $u$ is adjacent to $v$ in $\Lambda$. Therefore, the metric representation of the vertices $u, v\in V(\Lambda)-R$  is the same as $n$-vector, relative to $R$. Thus, $R$ is not a resolving set of $\Lambda$.
\newline

Case 3.
Now, let  $|R_1|=|R_2|$ and  suppose that for all the vertices  $x, y$ in $R_1$, we have $x$ is not adjacent to  $y$ in $\Lambda$, that is $d(x, y)=2$. Also for  all the vertices $u, v$ in $R_2$, we have $u$ is not adjacent to $v$ in $\Lambda$, that is $d(u, v)=2$. We may assume  that
$R_1=\{a,a^2, ..., a^{\frac{n}{2}}\}$ and  $R_2=\{ab, a^2b, ..., a^{\frac{n}{2}}b\}$.
So, we can assume that an arranged subset $R$ of vertices in graph $\Lambda$ is
$R=\{a, a^2, ..., a^{\frac{n}{2}}; ab, a^2b, ..., a^{\frac{n}{2}}b \}$.
Hence,
$V(\Lambda)- R=\{a^{\frac{n+2i}{2}}, ..., a^n; a^{\frac{n+2i}{2}}b, ..., a^nb\}$, for $1\leq i\leq \frac{n}{2}$.
Therefore, the metric representations of the vertices
$a^{\frac{n+2}{2}}, a^{\frac{n+4}{2}}, ..., a^n; a^{\frac{n+2}{2}}b, a^{\frac{n+4}{2}}b,..., a^nb  \in V(\Lambda)-R$
relative to $R$ are the $n$-vectors\newline
$r(a^{\frac{n+2}{2}} | R)=(1,2,2, ...,2; 1, ..., 1)$, $r(a^{\frac{n+4}{2}} |  R)=(2,1,2, ...,2; 1, ..., 1)$,  ..., $r(a^n |  R)=(2, 2,...,1 ; 1, 1,...1)$, and $r(a^{\frac{n+2}{2}}b |  R)=( 1, ..., 1; 1,2,2, ...,2 )$, $r(a^{\frac{n+4}{2}}b |  R)=( 1, ..., 1; 2,1,2, ...,2)$, ...,  $r(a^nb |  R)=( 1, 1,...1; 2, 2,...,1 )$. Thus, all the vertices in $V(\Lambda)-R$ have different representations relative to $R$. This implies that $R$ is a resolving set of $\Lambda$. \newline
\end{proof}
\begin{theorem}\label{d.2}
If $n$ is an even integer greater than or equal $4$, and $\Lambda=Cay(\mathbb{D}_{2n}, \Psi)$ be a Cayley graph on the dihedral group $\mathbb{D}_{2n}$, where $\Psi$  which is defined already, then  the cardinality of minimum doubly resolving set of $\Lambda$ is $n$.
\end{theorem}
\begin{proof}
By the previous Theorem, we know that the arranged subset
$R=\{a, a^2, ..., a^{\frac{n}{2}}; ab, a^2b, ..., a^{\frac{n}{2}}b\}$
of vertices in the graph $\Lambda$ is a resolving set for $\Lambda$. We show that the subset $R$ is a doubly resolving set of $\Lambda$. It is sufficient to show that for two  vertices $u$ and $v$ in graph $\Lambda$ there are vertices $x, y \in R$ such that $d(u, x) - d(u, y) \neq d(v, x) - d(v, y) $.
Consider two vertices $u$ and $v$ of $\Lambda$. By the following cases, we can be concluded  that the minimum cardinality of a doubly resolving set of $\Lambda$ is $n$.\newline

Case 1.
Consider a pair of distinct vertices $u, v\in \Lambda$ such that  $u, v\in R$. Then the length of a shortest path from $u$ to $v$  is $d(u, v)=1$ or $2$.
Let $u, v$ be two vertices in  $R$ such that a shortest path from $u$ to $v$ in graph $\Lambda$, is $d(u, v)=1$.
We may assume  that  $u=a$ and $v=ab$. Hence, by taking  $x=a\in R$ and $y=a^{\frac{n}{2}}\in R$, we have
$-2=0-2=d(u, x) - d(u, y) \neq d(v, x) - d(v, y)=1-1=0$.
Therefore, the vertices $x$ and $y$ of $R$ doubly resolve $u, v$. Now, let $u, v$ be two vertices in $R$ such that  a shortest path from $u$ to $v$ in  graph $\Lambda$, is $d(u, v)=2$.
We may assume  that  $u=a$ and $v=a^{\frac{n}{2}}$. Hence, by taking $x=a\in R$ and $y=ab\in R$ we have
$-1=0-1=d(u, x) - d(u, y) \neq d(v, x) - d(v, y)=2-1=1$.
Therefore, the vertices $x$ and $y$ of $R$ doubly resolve $u, v$.\newline

Case 2.
Consider a pair of distinct vertices $u, v\in \Lambda$ such that $u\in R$ and $v\notin R$. Then the length of a shortest path from $u$ to $v$  is $d(u, v)=1$ or $2$. Suppose a pair of distinct vertices $u\in R$ and $v\notin R$ are adjacent in  graph $\Lambda$ that is $d(u, v)=1$. We may assume  that $u=a$ and $v=a^{\frac{n+2}{2}}$. Hence,
by taking  $x=a\in R$ and $y=ab\in R$, we have $-1=0-1=d(u, x) - d(u, y) \neq d(v, x) - d(v, y)=1-1=0$.
Therefore, the vertices $x$ and $y$ of $R$ doubly resolve $u, v$. Now, suppose a pair of distinct vertices $u\in R$ and $v\notin R$ are not adjacent in  graph $\Lambda$, that is $d(u, v)=2$. We may assume  that $u=a$ and $v=a^n$. Hence, by taking $x=a\in R$ and $y=ab\in R$ we have $-1=0-1=d(u, x) - d(u, y) \neq d(v, x) - d(v, y)=2-1=1$.
Therefore, the vertices $x$ and $y$ of $R$ doubly resolve $u, v$.\newline

Case 3.
Consider a pair of distinct vertices $u, v\in \Lambda$ such that $u\notin R$ and $v\notin R$. Then the length of a shortest path from $u$ to $v$  is $d(u, v)=1$ or $2$. We can show that the subset $R$ of vertices in  graph $\Lambda$ is  a doubly resolving set of $\Lambda$. Because by Theorem 3.1, 
 we can be concluded that $V(\Lambda)-R$ is also resolving set of $\Lambda.$
\end{proof}
\begin{lemma}\label{d.3}
If $n$ is an even integer greater than or equal $4$, and $\Lambda=Cay(\mathbb{D}_{2n}, \Psi)$ be a Cayley graph on the dihedral group $\mathbb{D}_{2n}$, where $\Psi$  which is defined already, then the subset
$R=\{a, a^2, ..., a^{\frac{n}{2}}; ab, a^2b, ..., a^{\frac{n}{2}}b\}$ of vertices in  graph $\Lambda$ is not a strong resolving set of $\Lambda$.
\end{lemma}
\begin{proof}
We know that the arranged subset $R=\{a, a^2, ..., a^{\frac{n}{2}}; ab, a^2b, ..., a^{\frac{n}{2}}b\}$ of vertices in  graph $\Lambda$ is a resolving set for $\Lambda$ of size $n$. Now, let $V(\Lambda)=V_1 \cup V_2$, where $V_1=\{a,a^2, ..., a^n\}$, $V_2=\{ab,a^2b, ..., a^nb\}$, and $R=R_1\cup R_2$, where
$R_1=\{a,a^2, ..., a^{\frac{n}{2}}\}$ is a subset of $V_1$ and  $R_2=\{ab, a^2b, ..., a^{\frac{n}{2}}b\}$ is a subset of $V_2$. Consider two vertices $u, v$ in $\Lambda$ such that $u, v \in V_1-R_1$ and $u$ is not adjacent to  $v$ in $\Lambda$, that is $d(u, v)=2$.
In the following cases, we show  that there is not $w\in R$ such that $w$ is strongly resolves vertices $u$ and $v$. For every vertex $w\in R$, we have $w\in R_1$ or $w\in R_2$.
\newline

Case 1.
If $w\in R_1$, then the length of a shortest path from $u$ to $w$ is $d(u, w)=1$ or $2$, and length of a shortest path from $v$ to $w$ is $d(v, w)=1$ or $2$. Note that, if $d(u, w)=1$ then $d(v, w)=2$. Therefore, $w$ is not strongly resolves vertices $u$ and $v$. In particular, if $d(u, w)=2$ then $d(v, w)=1$ or $2$, and hence  $w$ is not strongly resolves vertices $u$ and $v$, because $d(u, v)=2$.\newline

Case 2.
If  $w\in R_2$, then the length of a shortest path from $u$ to $w$ is $d(u, w)=1$ and length of a shortest path from $v$ to $w$ is $d(v, w)=1$. Therefore, $w$ is not strongly resolves vertices $u$ and $v$.
\end{proof}
\begin{theorem}\label{d.4}
If $n$ is an even integer greater than or equal $4$, and $\Lambda=Cay(\mathbb{D}_{2n}, \Psi)$ be a Cayley graph on the dihedral group $\mathbb{D}_{2n}$, where $\Psi$  which is defined already, then the strong metric dimension of $\Lambda$ is $2n-2$.
\begin{proof}
Let $V(\Lambda)=V_1 \cup V_2$, where $V_1=\{a,a^2, ..., a^n\}$ and $V_2=\{ab,a^2b, ..., a^nb\}$.
It is not hard to see that if  $n\geq 4$, then the size of  largest clique in the graph $\Lambda$ is  $4$. Moreover, we know that the subset $N=\{a^n,  a^{\frac{n}{2}}; a^nb, a^{\frac{n}{2}}b\}$ of vertices in $\Lambda$ is a clique in the graph $\Lambda$. Now, let the subset $S$ of vertices in  $\Lambda$ is $S=V(\Lambda)-N$. In the following cases, we show  that the subset  $S$ of vertices in  $\Lambda$ is not a strong resolving set of $\Lambda$.\newline

Case1.
Let $u=a^n$, $v=a^{\frac{n}{2}}$. We know that $d(u,v)=1$,  $u,v\in V_1$, and hence for every $w\in S$ such that $w\in V_2$ we have $d(u, w)=1$ and $d(v, w)=1$.   Thus, $w$ is not strongly resolves vertices $u$ and $v$.\newline

Case 2.
Now, let $u=a^n$, $v=a^{\frac{n}{2}}$. We know that $d(u,v)=1$, $u,v\in V_1$, and hence for every $w\in S$ such that $w\in V_1$ we have $d(u, w)=2$ and $d(v, w)=2$.  Thus, $w$ is not strongly resolves vertices $u$ and $v$. \newline

Therefore, the subset  $S$ of the vertices in  graph $\Lambda$ is not  a strong resolving set of $\Lambda$. From  the above cases, we can be concluded that the minimum cardinality of a strong resolving set for $\Lambda$ must be $2n-2$.
\end{proof}
\end{theorem}
\section{Conclusion}
Computing resolving parameters of a graph is an NP-hard problem.  In this article, we considered a class of Toeplitz graphs, and we denoted by $T_{2n}(W)$, so that they are isomorphic to the Cayley graph
$\Lambda=Cay(\mathbb{D}_{2n}, \Psi)$, which is defined already. In fact, this class of graphs are vertex transitive,
and by calculating the spectrum of the adjacency matrix related with them, we showed that this class of  graphs cannot be edge transitive. Also, we proved  that this class of  graphs cannot be distance regular, and since the computing resolving parameters of a class of graphs such that are not distance regular is more difficult, then we regarded this as justification for our focus on some resolving parameters. In particular, we determined the minimal resolving set, doubly resolving set and strong metric dimension for this class of graphs.\newline

\bigskip
{\footnotesize
\noindent \textbf{Data Availability}\\
No data were used to support this study.\\[2mm]
\noindent \textbf{Conflicts of Interest}\\
The authors declare that there are no conflicts of interest
regarding the publication of this paper.\\[2mm]
\noindent \textbf{Acknowledgements}\\
This work was supported in part by Natural Science Fund of Education Department of Anhui Province under Grant KJ2020A0478.
\\[2mm]
\noindent \textbf{Authors' informations}\\
\noindent Jia-Bao Liu${}^a$
(\url{liujiabaoad@163.com;liujiabao@ahjzu.edu.cn})\\
Ali Zafari${}^{b}$(\textsc{Corresponding Author})
(\url{zafari.math.pu@gmail.com}; \url{zafari.math@pnu.ac.ir})\\

\noindent ${}^{a}$ School of Mathematics and Physics, Anhui Jianzhu University, Hefei 230601, P.R. China.\\
${}^{b}$ Department of Mathematics, Faculty of Science,
Payame Noor University, P.O. Box 19395-4697, Tehran, Iran.
}
\bigskip
{\footnotesize

\bigskip
\end{document}